\documentclass{article}
\pdfoutput=1

\usepackage{amsmath,amssymb,amsfonts}

\usepackage{amsthm,amsxtra}

\usepackage{graphicx}
\usepackage[svgnames]{xcolor}
\usepackage{verbatim}

\usepackage{pgf,tikz}
\usetikzlibrary{arrows.meta,decorations.pathmorphing,backgrounds,positioning,fit,petri,math}
\usepackage{tikz-bagua}

\usepackage{CJK}

\newcommand{\ggl}{\mathrm{GL}}
\newcommand{\ssl}{\mathrm{SL}}

\newcommand{\skxx}[1]{^{\mbox{\tiny{#1}}}}

\newcommand{\Vyo}{V\skxx{一}}

\newcommand{\Vyt}{V\skxx{二}}

\newcommand{\Vys}{V\skxx{三}}

\newcommand{\Vsh}{V\skxx{生}}
\newcommand{\Vke}{V\skxx{克}}

\newcommand{\Vsk}{V\skxx{生克}}

\newcommand{\Vtj}{V\skxx{太极}}

\newcommand{\Gpf}{G\skxx{平凡}}

\newcommand{\Gwx}{G\skxx{五行}}

\newcommand{\Gyy}{G\skxx{阴阳}}
\newcommand{\Gsx}{G\skxx{四象}}
\newcommand{\Gbg}{G\skxx{八卦}}

\newcommand{\Vwx}{V\skxx{五行}}
\newcommand{\Vly}{V\skxx{两仪}}

\newcommand{\Vsx}{V\skxx{四象}}
\newcommand{\Vbg}{V\skxx{八卦}}

\newcommand{\txt}[1]{{#1}}

\begin{document}

\title{Drawing the diagrams for Yinyang, Wuxing and Bagua as McKay quivers}
\author{Jin Yun Guo\\ School of mathematics and Statistics\\ Hunan Normal University}
\date{}

\maketitle

\begin{abstract}In this paper, we draw diagrams of yinyang, wuxing and bagua as the McKay quivers of the groups of order 2, of order $5$ and of order $2^3$.
\end{abstract}

\bigskip

\begin{CJK*}{UTF8}{gbsn} 

In \cite{m79}, McKay established the correspondence among the  finite subgroups of $\ssl(2,\mathbb C)$,  the simple complex Li algebras and Klein singularities. Now McKay correspondence is a much studied field in mathematics.
The finite subgroups of $\ssl(2,\mathbb C)$ are related to the finite subgroup of the rotation group $\mathrm{SO}(3,\mathbb R)$ of real space of dimension $3$. The $3$-dimensional geometric objects invariant under the finite subgroup of $\mathrm{SO}(3,\mathbb R)$ are exactly the Platoic solids.
McKay conclude his paper \cite{m79} with the words "Would not the Greeks appreciate the result that the simple Lie algebra may be derived from the Platonic solids?"

Yinyang(阴阳)\cite{r20}, wuxing(五行)\cite{f19} and bagua(八卦)\cite{h19} are conceptual frameworks in ancient Chinese thoughts in oberserving and analyzing the world. 
They are used even today in many aspects such as practising tranditional Chinese medicine, Fengshui, and so on.
There are diagrams illustrating the meaning for yinyang and wuxing.

Returning arrow quiver and covering of quiver are two constructions we used in constructing algebras in higher representation theory \cite{g13,gyz14,gl16}. We recover Taijitu for yinyang and the diagrams for wuxing using the McKay quivers for the groups of order $2$ and of order $5$, repsectiely, we also draw a  diagram for bagua using the McKay quiver of order $2^3$. 
Our approach gives interpretations to some assertions in {\em Daodejing} and in {\em Yijing}, from the point of view of representation theory.

Quivers are important realizations of certain categories.
By expressing the diagrams of yinyang, wuxing and bagua as quivers, we may rediscover the mathematical idea hide behind the old Chinese thoughts, under the point of view of category.

\section{Preliminary.}

In 1979, John McKay introduced McKay quiver in \cite{m79} for finite groups and their representations, and discovered the relationship between the Dynkin diagrams and the McKay quivers of finite subgroups of $\ssl(2,\mathbb C)$.

Let $k$ be an algebraically closed field of characteristics zero and let $V$ be a vector space of dimension $n$ over $k$.
Consider a finite  subgroup $G$ of the general linear group $\ggl(n,k)=\ggl(V)$, $V$ is naturally a representation of $G$.
Assume that the irreducible representations of $G$ are $S_1,\ldots,S_m$, and assume that we have a decomposition
$$V\otimes S_i = \bigoplus_{j=1}^m a_{i,j}S_j, $$
of the tensor product $V\otimes S_i$ for $ i=1,2,\cdots,m$.
{\em The McKay quiver} $Q_{V}(G)$ of $G$ is defined as the quiver with the vertices $S_1,\cdots,S_m$, and there are $a_{i,j}$ arrows from $S_i$ to $S_j$. 

In 2011, we prove the following theorems illustration the two constructions of McKay quivers in \cite{g11}:

{\bf Theorem 1.} (Theorem 1.2 in \cite{g11}) {\em Assume that $G$ is a finite subgroup of $\ggl(n,k)$ and that $N=G\cap \ssl(n,k)$. If every irreducible representation of $N$ is extendible, then the McKay quiver of $G$ is a regular covering of the McKay quiver of $N$ with $G/N$ as the group of the automorphisms.
}

The groups we considered are abelian, so the  irreducible representations of its subgroups are extendible.

\medskip

{\bf Theorem 2.} (Theorem 3.1 of \cite{g11}) {\em 
Assume that $G$ is a finite subgroup of $\ggl(n,k)$. 
Embed $\ggl(n,k)$ naturally into $\ssl(n+1,k)$.
Then the McKay quiver of $G$ in $\ssl(n+1,k)$ is obtained as the returning arrow quiver of the  McKay quiver of $G$ in $\ggl(n,k)$.} 

\medskip

Regarding the McKay quiver $Q_V(G)$ of $G\subset \ggl(V)$ as the bound quivers for the skew group algebra of $G$  over exterior algebra of $V$ \cite{gm02}, then non-zero paths in the algebra has finite length. For each maximal nonzero path $p$ in the algebra, adding an arrow, called "returing arrow", from the ending vertex of $p$ to the starting vertex of $p$. In this way, we get {\em the returning arrow quiver} of $Q_V(G)$ \cite{g11,gyz14}.

\medskip

In this paper, we first draw the diagrams of generating and of overcoming for wuxing as the McKay quivers of cyclic group of order $5$, directly from definition. 
We also show how to  construct the McKay quivers of  groups of order $2$ and of order $2^3$, using Theorems 1 and 2,  to get the diagram (Taijitu) for yinyang and a diagram for bagua， respectively. 


\section{Wuxing.}

Wuxing ("five processes" or "five phases", also translated as "five elements") refers to a fivefold conceptual scheme that is found throughout traditional Chinese thought.
These five phases are wood (mu), fire (huo), earth (tu), metal (jin), and water (shui); they are regarded as dynamic, interdependent modes or aspects of the universe’s ongoing existence and development \cite{wx}.  
The origins of wuxing trace back to Shang dynasty(1600-1046 B.C.E.), and was developed to its present form in Han dynasty((206 B. C. E.–220 C. E.).
The doctrine of wuxing is discribed in two diagrams, a diagram of generating or creation (生),  and a diagram of overcoming or destruction (克). We now draw the McKay quivers of the cyclic group of order $5$ as the  diagram of generating and the  diagram of overcoming of wuxing. 

\begin{figure}[htb]
\center
\begin{tikzpicture}[scale=0.5]

\def \radius {3cm}
\def \margin {8} 

  \node[draw, circle] (shui) at (18:\radius) {\txt{水}};

  \node[draw, circle] (jing) at (90:\radius) {\txt{金}}
  edge[post] node[midway]{\tiny\txt{生}} (shui);

  \node[draw, circle] (tu) at (162:\radius) {\txt{土}}
  edge[post] node[midway]{\tiny\txt{生}} (jing);

  \node[draw, circle] (huo) at (234:\radius) {\txt{火}}
  edge[post] node[midway]{\tiny\txt{生}} (tu);

  \node[draw, circle] (mu) at (305:\radius) {\txt{木}}
  edge[pre] node[midway]{\tiny\txt{生}} (shui)
  edge[post] node[midway]{\tiny\txt{生}} (huo);

\end{tikzpicture}
\caption{Wuxing--the diagram of generating
}\label{wu:sheng}
\end{figure}
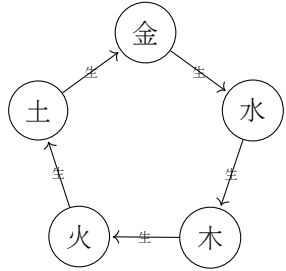

The cyclic group $\Gwx$ of order $5$ is generated by an element $a$ of order $5$. It has $5$ irreducible representations, all of dimension $1$. 
The values of the characters of these representations are the fifth roots of the unit.
The irreducible representations are determined by the value of their characters on the generator $a$ of $\Gwx$.
Denote by  $S_t$ the irreducible representation whose character takes value  $e^{\frac{2t\pi i}{5}}$ on $a$ for $0\le t\le 4$.
Write“金”$=S_{0}$,“水”$=S_{1}$,“木”$=S_{2}$,“火”$=S_{3}$,“土”$=S_{4}$.

Take $\Vwx=S_1$ and regard $\Gwx$ as a subgroup of $\ggl(1,k)=\ggl(\Vwx)$. Write “生” for tensoring with $ \Vsh=S_1$. 
We get the McKay quiver $Q_{\Vsh}(\Gwx)$(see Figure \ref{wu:sheng}) by direct calculation, which is the diagram of generating of wuxing.

\begin{figure}[htb]
\center
\begin{tikzpicture}[scale=0.5]

\def \radius {3cm}
\def \margin {8} 

  \node[draw, circle] (shui) at (18:\radius) {\txt{水}};

  \node[draw, circle] (jing) at (90:\radius) {\txt{金}};

  \node[draw, circle] (tu) at (162:\radius) {\txt{土}}
  edge[post] node[midway]{\tiny\txt{克}} (shui);

  \node[draw, circle] (huo) at (234:\radius) {\txt{火}}
  edge[post] node[midway]{\tiny\txt{克}} (jing)
  edge[pre] node[midway]{\tiny\txt{克}} (shui);

  \node[draw, circle] (mu) at (305:\radius) {\txt{木}}
  edge[post] node[midway]{\tiny\txt{克}} (tu)
  edge[pre] node[midway]{\tiny\txt{克}} (jing);

\end{tikzpicture}
\caption{Wuxing--the diagram of overcoming 
}\label{wu:ke}
\end{figure}
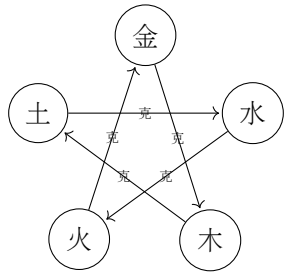

Now take $\Vwx=S_2$, write “克” for tensoring with $\Vke=S_2$.
We get the McKay quiver $Q_{\Vke}(\Gwx)$ (see Figure \ref{wu:ke}), this  is the diagram of overcoming of wuxing  when Figure \ref{wu:sheng}  is the diagram of generating of wuxing.

When $\Vwx$ is taken as $S_{4}$ or $S_3$, the McKay quivers become the diagrams of weakening (泻，inverse of generating) and the diagram of insulting (侮，inverse of overcoming).

\begin{figure}[htb]
\center
\begin{tikzpicture}[scale=0.5]

\def \radius {3cm}
\def \margin {8} 

  \node[draw, circle] (shui) at (18:\radius) {\txt{水}};

  \node[draw, circle] (jing) at (90:\radius) {\txt{金}}
  edge[post] node[midway]{\tiny\txt{生}} (shui);

  \node[draw, circle] (tu) at (162:\radius) {\txt{土}}
  edge[post] node[midway]{\tiny\txt{生}} (jing)
  edge[post] node[midway]{\tiny\txt{克}} (shui);

  \node[draw, circle] (huo) at (234:\radius) {\txt{火}}
  edge[post] node[midway]{\tiny\txt{克}} (jing)
  edge[pre] node[midway]{\tiny\txt{克}} (shui)
  edge[post] node[midway]{\tiny\txt{生}} (tu);

  \node[draw, circle] (mu) at (305:\radius) {\txt{木}}
  edge[post] node[midway]{\tiny\txt{克}} (tu)
  edge[pre] node[midway]{\tiny\txt{克}} (jing)
  edge[pre] node[midway]{\tiny\txt{生}} (shui)
  edge[post] node[midway]{\tiny\txt{生}} (huo);

\end{tikzpicture}
\caption{Wuxing--the diagram of generating and overcoming 
}\label{wu:shengke}
\end{figure}
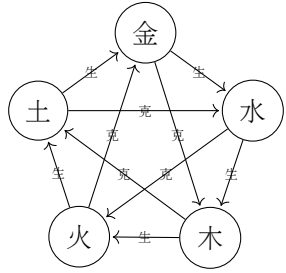

Now take $\Vsk=S_1\oplus S_2$ and regard $\Gwx$ as a subgroup of $\ggl(2,k)=\ggl(\Vsk)$,  we get the McKay quiver $Q_{\Vsk}(\Gwx)$（Figure \ref{wu:shengke}), which is the generating-overcoming diagram of wuxing.

So the diagram of generating and the diagram of overcoming of wuxing are exactly the McKay quivers of cyclic group of order $5$.

\section{Yinyang.}

Yinyang (yin-yang) is one of the dominant concepts shared by different schools throughout the history of Chinese thoughts.
The Taijitu (Diagram of the Great Ultimate) for  yinyang attributed to Zhou Dunyi (1017-1073 CE)\cite{yy}. 
Following the guidance of Chapter 42 of {\em Daodejin} (\cite{ddj}), we draw the diagram of {\em Taijitu} as McKay quivers, applying Theorems 1 and 2.

\begin{figure}[htb]
\center
\begin{tikzpicture}[scale=1.4]
\filldraw [green] (0,0) circle (3pt);
\draw[gray, ultra thick, ->] (0,0) arc (90:445:1);
\end{tikzpicture}
\caption{Taiji}\label{yy:taiji}
\end{figure}

The simplest group is the trivial group with only one element, the unit.
The next simple group is the group of order $2$.
It has two elements, the one other than the unit is a generator.
Any group has a trivial irreducible representation $S_0$. It is one dimensional, and the value of its character on any element is $1=e^0$.

Take $\Vtj=S_0$, and we consider the finite subgroup of $\ssl(1,k)$. The only such group is the trivial group $\Gpf$, the group with only one element. Now we draw its McKay quiver $Q_{\Vtj}(\Gpf)$(Figure \ref{yy:taiji}). It is the quiver with only one vertex and a loop arrow on the vertex, since $S_0\otimes S_0= S_0$. We view as this quiver as Taiji(太极，the supreme ultimate).

\begin{figure}[htb]
\center
\begin{tikzpicture}[scale=1.4]

\filldraw [green] (0,1) circle (3pt) node[anchor=south]{\txt{阳}};
\filldraw [green] (0,-1) circle (3pt) node[anchor=north]{\txt{阴}};
\draw[gray, ultra thick, ->] (0,1) arc (90:270:1);
\draw[gray, ultra thick, ->] (0,-1) arc (270:450:1);
\end{tikzpicture}
\caption{Covering quiver
}\label{yy:yse}
\end{figure}
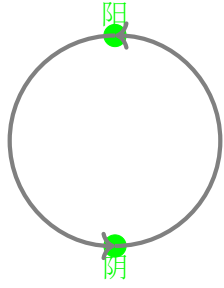

The trivial group only has the unit, which is written as $0$ in general, and means "has nothing", or "none".
As a group it has a trivial representation, that is, it "has something": "and being is born of nothing"(Chapter 40 of \cite{ddj}).

Write $\Gyy $ for the cyclic group of order $2$, and write $a$ for its generator.  It has two irreducible representations, $S_0$ and and $S_1$, whose characters  take value $1=e^0$ and $e^{\pi i}=-1$, respectively, at $a$. Write 阳(Yang)=$S_1$ and 阴(Yin)=$S_0$ for these two irreducible representations.

Take $\Vyo=S_1$, then $\Gyy$ is a finite subgroup of $\ggl(1,k)=\ggl(\Vyo)$. In this case, $\Gyy\cap \ssl(1,k)= \Gpf$, and $\Gyy/\Gpf \simeq  \Gyy$. By Theorem 1, the McKay quiver $Q_{\Vyo}(\Gyy)$(Figure \ref{yy:yse}) of $\Gyy$ is a regular cover of the McKay quiver of $\Gpf$ with $\Gyy$ as the group of automorphisms. Since $G$ has two elements, under its action, we get two vertices 阳 and 阴 from the vertex in $Q_{\Vtj}(\Gpf)$ and two arrows in the opposition direction between the two vertices from the loop arrow in   $Q_{\Vtj}(\Gpf)$.

According to \cite{gm02,gxl21}, we can associate to this McKay quiver a $0$-slice algebra which is of finite representation type.

We remark that Figure \ref{wu:sheng} and  Figure \ref{wu:ke} for wuxing can also be obtained using this approach.

\begin{figure}[htb]
\center
\begin{tikzpicture}[scale=0.7]

\filldraw [green] (0,2) circle (5pt) node[anchor=south]{\txt{阳}};
\filldraw [green] (0,-2) circle (5pt) node[anchor=north]{\txt{阴}};
\draw[gray, ultra thick, ->] (0,2) arc (90:270:2);
\draw[gray, ultra thick, ->] (0,-2) arc (270:450:2);
\draw[red, thick, -{Latex[right]}] (0,2) .. controls(-0.5,0) .. (0,-1.9) ;
\draw[red, thick, -{Latex[right]}] (0,-2) .. controls(0.5,0) .. (0,1.9) ;

\end{tikzpicture}
\caption{Returning arrows--``One gives birth to two''
}\label{yy:htj}
\end{figure}
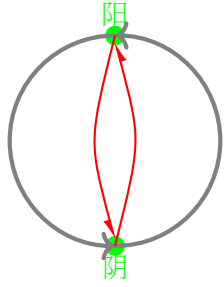

Take $\Vyt=S_1\oplus S_1$ and we have a vector space of dimension $2$. The group $\ggl(1,k)$ is embedded into $\ssl(2,k)=\ssl(\Vyt)$. With this  embedding, $\Gyy$  is now a subgroup of $\ssl(2,k)$.
By Theorem 2, the McKay quiver $Q_{\Vyt}(\Gyy)$ (see  Figure \ref{yy:htj}) of $\Gyy$ in $\ssl(\Vyt)$ is obtained from  the McKay quiver $Q_{\Vyo}(\Gyy)$ of $\Gyy$ in $\ggl(\Vyo)$ by adding at each vertex a "returning arrow". It follows from \cite{gm02,g11} that the "returning arrow" starting at a vertex is just an arrow opposite to the arrow ending at this vertex in $Q_{\Vyo}(\Gyy)$. Now we have two arrows from each vertex, in this way, we see that "one gives birth to two".

It follows from \cite{gm02,gxl21} that the $1$-slice algebra associated to this McKay quiver is the Kronecker algebra, which is  of tame representation type, meaning that its indecomposable representations can be parametrized   by the projective space of dimension $1$.

\begin{figure}[htb]
\center
\begin{tikzpicture}[scale=0.7]

\filldraw [green] (0,2) circle (5pt) node[anchor=south]{\txt{阳}};
\filldraw [green] (0,-2) circle (5pt)node[anchor=north]{\txt{阴}};
\draw[gray, ultra thick, ->] (0,2) arc (90:270:2);
\draw[gray, ultra thick, ->] (0,-2) arc (270:450:2);
\draw[purple!50, thick, -{Stealth[right]}] (0,2) .. controls(-0.5,0) .. (0,-1.9) ;
\draw[purple!50, thick, -{Stealth[right]}] (0,-2) .. controls(0.5,0) .. (0,1.9) ;
\draw[blue, ultra thick, -Latex] (0,2) arc (90:450:0.5);
\draw[blue, ultra thick, -Latex] (0,-2) arc (270:630:0.5);

\end{tikzpicture}
\caption{"Returning arrows" in dimension $3$--"two gives birth to three''
}\label{yy:swhtj}
\end{figure}
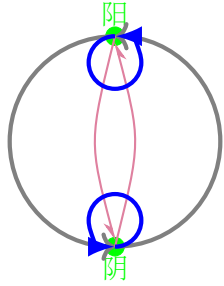

Now take $\Vys=S_0\oplus S_1\oplus S_1$, we get a vector space of dimension $3$. In this case, the group $\ssl(2,k)$ is embedded into $\ssl(3,k) =\ssl(\Vys) $ and $\Gyy$ becomes a finite subgroup of $\ssl(3,k)$. Again by Theorem 2, the McKay quiver $Q_{\Vys}(\Gyy)$ of $\Gyy$ in  $\ssl(3,k)$ is obtained from the McKay quiver $Q_{\Vyt}(\Gyy)$ by adding at each vertex a "returning arrow". In this case, the returning arrows are loops  at every vertices, since we have added a trivial irreducible representation to the space (see \cite{g11}). So we get the McKay quiver of $\Gyy$ in dimension $3$ as $Q_{\Vys}(\Gyy)$(Figure \ref{yy:swhtj}).

\begin{figure}[htb]
\center
\begin{tikzpicture}[scale=0.7]
\filldraw [green] (0,2) circle (5pt)  node[anchor=south]{\txt{阳}};
\filldraw [green] (0,-2) circle (5pt)  node[anchor=north]{\txt{阴}};
\draw[gray, ultra thick, ->] (0,2) arc (90:270:2);
\draw[gray, ultra thick, ->] (0,-2) arc (270:450:2);

\draw[blue, ultra thick, -Latex] (0,2) arc (90:450:0.5);
\draw[blue, ultra thick, -Latex] (0,-2) arc (270:630:0.5);
\draw[purple,  thick, double distance=1pt, -{Stealth[right]}
] (0,0) arc (270:450:1);
\draw[purple,  thick, double distance=1pt, -{Stealth[right]}
] (0,0) arc (90:270:1);

\end{tikzpicture}
\caption{"Returning arrows" in dimension $3$ without crossing 
}\label{yy:swhtjwj}
\end{figure}
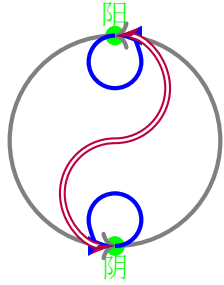

Each vertex of this McKay quiver is the starting vertex of $3$ arrows and also is the ending  vertex of $3$ arrows.
When the dimension rised to $3$ from $2$, two arrows give birth to three arrows, "two gives birth to three".

\begin{figure}[htb]
\center
\begin{tikzpicture}[scale=0.7]
 \begin{scope}
    \clip (0,0) circle (2cm);
    \fill[white]  (0,0) circle (2cm);
    \fill[black] (0cm,2cm) rectangle (-2cm, -2cm);
  \end{scope}

  \fill[black] (0,1) circle (1cm);
  \fill[white] (0,-1) circle (1cm);

  \fill[white] (0,1.5) circle (0.5cm);
  \fill[black] (0,-1.5) circle (0.5cm);

\filldraw [green!50] (0,2) circle (5pt)  node[anchor=south]{\txt{阳}};
\filldraw [green!50] (0,-2) circle (5pt)  node[anchor=north]{\txt{阴}};
\draw[gray, ultra thick, ->] (0,2) arc (90:270:2);
\draw[gray, ultra thick, ->] (0,-2) arc (270:450:2);

\draw[purple!50, ultra thick, double distance=1pt, -{Stealth[right]}
] (0,0) arc (270:450:1);
\draw[purple!50, ultra thick, double distance=1pt, -{Stealth[right]}
] (0,0) arc (90:270:1);

\draw[blue, ultra thick, -Latex] (0,2) arc (90:450:0.5);
\draw[blue, ultra thick, -Latex] (0,-2) arc (270:630:0.5);

\end{tikzpicture}

\caption{Taijitu}\label{yy:yytjt}
\end{figure}

Adjust the returning arrows so that they do not crossing with other arrows, we get Figure \ref{yy:swhtjwj}。
By \cite{hjc13}, this McKay quiver lies on the surface of a torus.
Fill it with black and white, we get Figure \ref{yy:yytjt}, which is the diagram of {\em Taijitu}. 

Using \cite{gm02,gxl21}, we associate a $2$-slice algebra with this McKay quiver. This algebra is of wild representation type, meaning that the category of the representations of any finite dimensional algebra can be embedded into the category of the representations of this algebra. So "three gives birth to the ten thousand things".

Our process of drawing McKay quivers to get the diagram of Taijitu reproduces what  Chapter 40 of {\em Daodejing} says: "The things of the world are born from being, and being is born of nothing"(see \cite{ddj}), and what Chapter 42 of {\em Daodejing} says "The Dao gives birth to one; one gives birth to two; two gives birth to three; three gives birth to the ten thousand things".

\section{Bagua.}

Bagua， the 8 trigrams (bagua of $3$ yaos) and 64 hexagrams(bagua of 6 yaos) allegedly created by the mythical figure Fuxi as the first layer of the  {\em Yijing} ({\em I Ching}, is also translated as {\em (Zhou) Book of Change(s)} \cite{zy}), a  canonized Confucian classic. Gottfried Wilhelm Leibniz, the inventor of binary  system, thought he had discovered forgotten mathematical science when he saw these diagrams in the order of Shao Yong(see \cite{r96}).

In mathematics, especially in representation theory, groups appear as a closed set of the transformations of some spaces. "Yi" is translated as "change". If we take "yi" as transformation, or call a transformation a "yi". As a study of "yi", the "yi"s in the book {\em Yijing} should form a closed set. So as a whole, bagua should have certain property related to group. As we will show, they form McKay quivers of some group with higher dimensional forms of cyclic group of order $2$, -- direct sums of three $\Gyy$'s and of six $\Gyy$'s, respectively.

In the first part of {\em The Survey} of {\em The Appendices} of {\em Yijing}(《易传·系辞上传》), it says  "Therefore in The Zhou Book of Change, there is the supreme ultimate(Taiji, 太极) at first. The supreme ultimate gives birth to the two primary forms (yin and yang: 两仪, Liangyi), the two primary forms give birth to the four basic images (greater yin, greater yang, lesser yin, less yang: 四象, Sixiang), and four basic images to eight trigrams(八卦, bagua)"(see \cite{zy}).
Now we show how bagua is realized under such process as McKay quivers, using Theorems 1 and 2.

We have seen in Figure \ref{yy:taiji} that Taiji is obtained as the McKay quiver of the trivial group embedded in one dimensional special linear group as a subgroup.

Now start with the cyclic group  $\Gyy$ of order $2$, we will show how bagua is reached using McKay quiver construction.
 
Write the representations of the group $\Gyy$ in another way, write the representations 阴 and 阳 as \liangyi{0} and \liangyi{1}, respectively. Take $\Vly=\Vyo$，the McKay quiver $Q_{\Vyo}(\Gyy)$(Figure \ref{yy:yse}) now become "Liangyi" quiver $Q_{\Vly}(\Gyy)$ (Figure \ref{bg:coveringly}).
\begin{figure}[htb]
\center
\begin{tikzpicture}[scale=0.7]

\node [draw,gray, ultra thick, circle] at (0,2){\txt{\liangyi{1}}};
\node [draw,gray, ultra thick, circle] at (0,-2) {\txt{\liangyi{0}}};
\draw[gray, ultra thick, ->] (-0.3,2) arc (98:262:2);
\draw[gray, ultra thick, ->] (0.3,-2) arc (278:442:2);
\end{tikzpicture}
\caption{Liangyi}\label{bg:coveringly}
\end{figure}

In Cartesian coordinates, a point in the space is expressed by a ordered set of numbers, and spaces of higher dimension are direct sums of lower dimensional spaces.
Apply this idea on groups to define the direct sum (or direct product when using multiplication for the compositions) of groups.
Take direct sum of two cyclic groups $\Gyy$ of order $2$, we obtain an abelian group $\Gsx= C_2\oplus C_2$, this is the Klein's $4$-group with $4$ elements $\{(0,0),(0,a),(a,0),(a,a)\}$.
$a_1=(a,0), a_2= (0,a)$ form a set of generators  of $\Gsx$. 
The values of the characters of the irreducible representations of $\Gsx$ are $\pm 1$, write \liangyi{0} when the value is $0$ and $\liangyi{1}$ when the value is $1$. Put the value of the character at $a_1$ in the position of first yao(first position from the bottom) and   the value of the character at $a_2$ in the position of the second yao(second position from the bottom).
The four irreducible representations of this group are now denoted by "Sixiang": \sixiang*{0}(太阴,great yin), \sixiang*{1}(少阳,lesser yang), \sixiang*{2}(少阴, lesser yin) and \sixiang*{3}(太阳, great yang).
    
\begin{figure}[htb]
\center
\begin{tikzpicture}[scale=0.7]

\node [draw,gray, ultra thick, circle] at (0,2) (ly1) {\txt{\liangyi{1}}};
\node [draw,gray, ultra thick, circle] at (0,-2) (ly2) {\txt{\liangyi{0}}} 
    edge[pre, purple!50, ultra thick, bend left=10]  (ly1)    
    edge[post, purple!50, ultra thick, bend right=10]  (ly1);
\draw[gray, ultra thick, ->] (-0.3,2) arc (98:262:2);
\draw[gray, ultra thick, ->] (0.3,-2) arc (278:442:2);
\end{tikzpicture}
\caption{Liangyi in 2D}\label{bg:lyin2}
\end{figure}
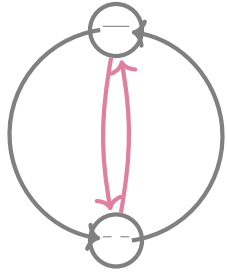

Take $\Vsx = \txt{\sixiang*{2}}\oplus \txt{\sixiang*{3}}$, we get a space of dimension $2$. 
Embed $\Gyy$ naturally into $\ssl(\Vsx)=\ssl(2,k)$ and get the McKay quiver $Q_{\Vsx}(\Gyy)$ of $\Gyy$ in $\ssl(2,k)$, a quiver with  "returning arrows" Figure \ref{bg:lyin2}, by Theorem 2.

Now $\Gsx$ is a subgroup of $\mathrm{GL}(2,\mathbb C)$ obtained by extending the subgroup $\Gyy$ of $\mathrm{SL}(2,\mathbb C)$ with a cyclic group $\Gyy$ of order $2$. By Theorem 1， its McKay quiver $Q_{\Vsx}(\Gsx)$ is a regular covering of $Q_{\Vsx}(\Gyy)$ with $\Gyy$ as the group of automorphisms. So under the action of $\Gyy$, the vertex \liangyi{0} becomes two vertices \sixiang*{0} and \sixiang*{2}, the vertex \liangyi{1} becomes two vertices \sixiang*{3} and \sixiang*{1}. The pair of old arrows becomes two pair of arrows in the opposite direction inside each of these pair of vertices, the pair of returning arrows becomes two pair of arrows in the opposite direction between the vertices from the different pairs. So we get Figure \ref{bg:sixiang} as the McKay quiver $Q_{\Vsx}(\Gsx)$.

\begin{figure}[htb]
\center
\begin{tikzpicture}[scale=0.7]

\def \radius {3cm}
\def \margin {8} 

  \node[draw, gray, ultra thick, circle] (sx0) at (225:\radius) {\txt{\sixiang*{0}}};

  \node[draw, gray, ultra thick, circle] (sx1) at (135:\radius) {\txt{\sixiang*{2}}}
    edge[post,gray, ultra thick, bend right=10]  (sx0)
    edge[pre, gray, ultra thick, bend left=10] (sx0);

  \node[draw, gray, ultra thick,  circle] (sx2) at (315:\radius) {\txt{\sixiang*{3}}}
    edge[post,purple!50, ultra thick, bend left=10] (sx0)    
    edge[pre, purple!50, ultra thick, bend right=10] (sx0);

  \node[draw, gray, ultra thick, circle] (sx3) at (45:\radius) {\txt{\sixiang*{1}}}
    edge[pre, purple!50, ultra thick, bend left=10]  (sx1)    
    edge[post, purple!50, ultra thick, bend right=10]  (sx1)
    edge[pre, gray, ultra thick, bend right=10] (sx2)    
    edge[post, gray, ultra thick, bend left=10] (sx2);
 
\end{tikzpicture}
\caption{Sixiang}\label{bg:sixiang}
\end{figure}

Add another direct summand $\Gyy$ on the Klein $4$-gruop $\Gsx$, we get bagua gruop $\Gbg$ with $8$ elements. 
$\Gbg$ needs $3$ generators, $a_1=(a,0,0), a_2=(0,a,0), a_3=(0,0,a)$. The group $\Gbg$ has $8$ irreducible representations. The characters of the  irreducible representations of $\Gbg$ take values $\pm 1$. 
The characters are determined by their values on the generators and we can use bagua to denote the values of the characters on the generators.  Put the value of the character at $a_1$ in the position of first yao (first position from the bottom),  the value of the character at $a_2$ in the position of the second yao (second position from the bottom), and   the value of the character at $a_3$ in the position of the third yao(the top).
The eight irreducible representations of the group $\Gbg$ are now denoted by "bagua": \bagua{000}(坤, kun)), \bagua{001}(艮, gen), \bagua{010}(坎, kan), \bagua{011}(巽, xun), \bagua{100}(震, zhen), \bagua{101}(离, li), \bagua{110}(兑, dui), \bagua{111}(乾, qian).

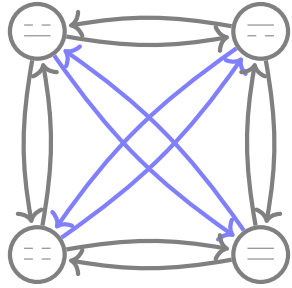
\begin{figure}[htb]
\center
\begin{tikzpicture}[scale=0.7]

\def \radius {3cm}
\def \margin {8} 

  \node[draw, gray, ultra thick, circle] (sx0) at (225:\radius) {\txt{\sixiang*{0}}};
 
  \node[draw, gray, ultra thick, circle] (sx1) at (135:\radius) {\txt{\sixiang*{2}}}
    edge[post,gray, ultra thick, bend right=10]  (sx0)
    edge[pre, gray, ultra thick, bend left=10] (sx0);

  \node[draw,  gray, ultra thick, circle] (sx2) at (315:\radius) {\txt{\sixiang*{3}}}
    edge[pre, blue!50, ultra thick, bend left=10]  (sx1)    
    edge[post, blue!50, ultra thick, bend right=10]  (sx1)
    edge[post,gray, ultra thick, bend left=10] (sx0)    
    edge[pre, gray, ultra thick, bend right=10] (sx0);

  \node[draw, gray, ultra thick, circle] (sx3) at (45:\radius) {\txt{\sixiang*{1}}}
    edge[pre, gray, ultra thick, bend left=10]  (sx1)    
    edge[post, gray, ultra thick, bend right=10]  (sx1)
    edge[pre, blue!50, ultra thick, bend left=10]  (sx0)    
    edge[post, blue!50, ultra thick, bend right=10]  (sx0)
    edge[pre, gray, ultra thick, bend right=10] (sx2)    
    edge[post, gray, ultra thick, bend left=10] (sx2);

\end{tikzpicture}
\caption{Sixiang: Returning arrow}\label{bg:sixiangr}
\end{figure}

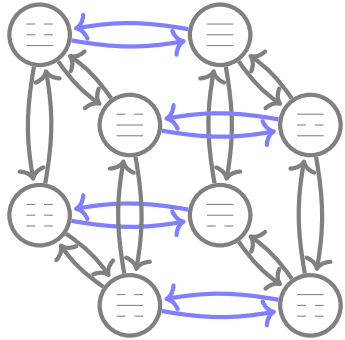
\begin{figure}[htb]
\center

\begin{tikzpicture}[scale=0.8]

\def \radius {3cm}
\def \margin {8} 

  \node[draw, gray, ultra thick, circle] (bg0) at (0,0) {\txt{\bagua{000}}};

  \node[draw, gray, ultra thick, circle] (bg1) at (3,0) {\txt{\bagua{011}}}
    edge[post,blue!50, ultra thick, bend right=10]  (bg0)
    edge[pre, blue!50, ultra thick, bend left=10] (bg0);

  \node[draw, gray, ultra thick, circle] (bg2) at (1.5,-1.5) {\txt{\bagua{010}}}
    edge[post, gray, ultra thick, bend left=10] (bg0)    
    edge[pre, gray, ultra thick, bend right=10] (bg0);
 
  \node[draw, gray, ultra thick, circle] (bg4) at (0,3) {\txt{\bagua{100}}}
    edge[pre, gray, ultra thick, bend left=10]  (bg0)    
    edge[post, gray, ultra thick, bend right=10]  (bg0);

  \node[draw, gray, ultra thick, circle] (bg3) at (4.5,-1.5) {\txt{\bagua{001}}}
    edge[post,gray, ultra thick, bend right=10]  (bg1)
    edge[pre, gray, ultra thick, bend left=10] (bg1)
    edge[post,blue!50, ultra thick, bend right=10]  (bg2)
    edge[pre, blue!50, ultra thick, bend left=10] (bg2);

  \node[draw, gray, ultra thick, circle] (bg5) at (3,3) {\txt{\bagua{111}}}
    edge[post, gray, ultra thick, bend left=10] (bg1)    
    edge[pre, gray, ultra thick, bend right=10] (bg1)
    edge[post,blue!50, ultra thick, bend right=10]  (bg4)
    edge[pre, blue!50, ultra thick, bend left=10] (bg4);

  \node[draw, gray, ultra thick, circle] (bg6) at (1.5,1.5) {\txt{\bagua{110}}}
    edge[post, gray, ultra thick, bend left=10] (bg2)    
    edge[pre, gray, ultra thick, bend right=10] (bg2)
    edge[post,gray, ultra thick, bend right=10]  (bg4)
    edge[pre, gray, ultra thick, bend left=10] (bg4);

  \node[draw, gray, ultra thick, circle] (bg7) at (4.5,1.5) {\txt{\bagua{101}}}
    edge[post, gray, ultra thick, bend left=10] (bg3)    
    edge[pre, gray, ultra thick, bend right=10] (bg3)
    edge[post,gray, ultra thick, bend right=10]  (bg5)
    edge[pre, gray, ultra thick, bend left=10] (bg5)      
    edge[pre, blue!50, ultra thick, bend left=10]  (bg6)    
    edge[post, blue!50, ultra thick, bend right=10]  (bg6);
\end{tikzpicture}

\caption{Bagua with $3$ yaos}\label{bg:sybg}
\end{figure}

Write $\Vbg=\txt{\bagua{100}}\oplus \txt{\bagua{110}}\oplus \txt{\bagua{011}} $, we ge a vector space of dimension $3$. It is obtained from $\Vsx$ by adding an one-dimensional subspace as a direct summand. Identify $\Vsx$ with the subspace $\txt{\bagua{100}}\oplus \txt{\bagua{110}}$, then $\Gsx$ is embedded into $\ssl(3,k)$. By Theorem 2, its McKay quiver (Figure \ref{bg:sixiangr}) is obtained from $Q_{\Vsx}(\Gsx)$ by adding a "returning arrow" at each vertex. In this case, such arrows are oppisite to the non-cyclic paths of length $2$. The group $\Gbg$ 
is a subgroup of $\ggl(3,k)$ obtained as the extension of $\Gsx$, embedded in $\ggl(2, k)$ as a subgroup in $\ssl(3,k)$, by a group $\Gyy$ of order $2$. So $\Gbg \cap \ssl(3,k) =\Gsx$ and $\Gbg/\Gsx \simeq \Gyy$.
By Theorem 1, $Q_{\Vbg}(\Gbg)$ is a regular covering of $Q_{\Vbg}(\Gsx)$ with $\Gyy$ as the group of automorphisms.
The action of $\Gyy$ on the vertices generates two sets \bagua{000}, \bagua{010}, \bagua{100}, \bagua{110} and \bagua{001}, \bagua{011}, \bagua{101}, \bagua{111} of vertices, it also generates two sets of arrows from the old arrows in $Q_{\Vsx}(\Gsx)$ such that they form two quiver $Q_{\Vsx}(\Gsx)$. The action of $\Gyy$ on the returning arrows produces arrows connecting the two quivers. In this way we get Figure \ref{bg:sybg}, the McKay quiver $Q_{\Vbg}(\Gbg)$ of $\Gbg$ in $\ggl(3,k)$ (See proof of Proposition 4.1 of \cite{gxl21} for a general construction).

This McKay quiver can be regarded as a unit cube in the space of dimension $3$, with each edge replaced by a pair of arrows heading in the opposite direction.

Similarly we can draw the quiver for the 64 hexagrams as McKay quiver for group which is direct sum of $6$ cyclic group $\Gyy$ of order $2$. 
It is too complecated, we will try to draw it nicely and give intepretation in the furture.

\bigskip

Now we have drawn the diagram of yinyang, wuxing and bagua as McKay quivers, may be representation theory can help to discover and develop the thought of the great mind in ancient time.

\medskip

{\bf Acknowledgement}. The author is supported by the  National Natural Science Foundation of China \#12071120.

{}

\end{CJK*}

\end{document}